\documentclass[10pt]{amsart} 

\usepackage{amsmath}
\usepackage{amssymb}
\usepackage{amsfonts}
\usepackage{graphicx}
\usepackage{epsfig}

\title{Voronoi Cells of Discrete Point Sets}

\author{Ina Voigt}
\address{I. Voigt, Faculty of Mathematics, TU Dortmund University, Vogelpothsweg 87, D-44227 Dortmund, Germany} 
\email{ina.voigt@mathematik.tu-dortmund.de}


\date{November 10, 2008}

\newtheorem{Def}{Definition}[section]
\newtheorem{Thm}{Theorem}[section]
\newtheorem{Cor}{Corollary}[section]
\newtheorem{Ex}{Example}[section]

\def\Hyp{\mathop{\mathrm{ Hyp}}\nolimits}      
\def\conv{\mathop{\mathrm{ conv}}\nolimits}

\def\aff{\mathop{\mathrm{ aff}}\nolimits}     

\def\Ext{\mathop{\mathrm{ Ext}}\nolimits}      
\def\RExt{\mathop{\mathrm{ RExt}}\nolimits}      
\def\cone{\mathop{\mathrm{ cone}}\nolimits}

\begin{document}

\begin{abstract}
It is well known that all cells of the Voronoi diagram of a Delaunay set are polytopes. For a finite point set, all these cells are still polyhedra. 
So the question arises, if this observation holds for all discrete point sets: Are always all Voronoi cells of an arbitrary, infinite discrete point set polyhedral?
In this paper, an answer to this question will be given.
It will be shown that all Voronoi cells of a discrete point set are polytopes if and only if every point of the point set is an inner point. Furthermore, the term of a locally finitely generated discrete point set will be introduced and it will be shown that exactly these sets have the property of possessing only polyhedral Voronoi cells.
\end{abstract}

\maketitle

\markboth{Ina Voigt}{Voronoi cells of discrete point sets}

\section{Introduction}
A review of the literature on Voronoi cells reveals the fact that they emerge in many different fields of research and are known under many different names like \emph{nearest neighbor region}, \emph{Wigner-Seitz cell} or \emph{Thiessen polygon}. 
In all these fields, the research activities on Voronoi diagrams were very vivid during the last decades.
A comprehensive overview is given e.g. in \cite{Auren00} or \cite{OBS00}, where a great part of the literature about Voronoi cells is cited.
Thereby one notices that in all applications, the related discrete point sets are either finite or possess a certain "structure".
The literature on computational geometry (see e.g. \cite{deBerg08} for an overview) is mainly concerned with the algorithmic construction of Voronoi diagrams. Hence, it is normally assumed that the set of considered point sites is finite.
On the other hand, if Voronoi diagrams of infinite point sets are observed, the examined set of points is always assumed to exhibit some sort of structure. 
In various fields of natural and social sciences one is interested in the so called \emph{Poisson Voronoi diagram} (cf. Chapter 5 of \cite{OBS00}), where the considered point set is a realization of a non-empty stationary point process and is usually assumed to be discrete and in \emph{general quadratic position}, which means that in $m$-dimensional Euclidean space no $k+1$ points lie on a $(k-1)$-dimensional hyperplane for $k = 2, \dots, m$ and that no $m+2$ points lie on the boundary of a sphere. 
In crystallography (see e.g. \cite{Baa02}), the considered point sets are generally assumed to be \emph{Delaunay sets} (to be defined in Section 2). 
It is well known that in all these cases all Voronoi cells are either polytopal or polyhedral. 
But what happens if we consider an \emph{arbitrary} infinite discrete point set? Are there still all Voronoi cells polyhedra or even polytopes?

In this paper we will show that this is not the case.
We will give an example of an infinite discrete point set whose Voronoi diagram contains a non-polyhedral cell, cf. Example \ref{Ex:non-polyhedral-cell}.
But obviously, there exist discrete point sets for which all Voronoi cells are polyhedral. 
Thus, the aim is to characterize the discrete point sets that possess only polyhedral or rather polytopal Voronoi cells.
By partitioning the discrete point set in \emph{inner points} and \emph{boundary points} relative to its convex hull (cf. Definition \ref{Def:inner+boundary-points}) and introducing the term of the \emph{direction cone} of a point of a discrete point set (cf. Definition \ref{C(p)}), we can formulate the desired characterizations:
In Theorem \ref{Thm1} we show that the Voronoi cell of a point $p$ is a polytope if and only if $p$ is an inner point of the discrete point set.
Theorems \ref{Thm2} and \ref{Thm3} prove that the Voronoi cell of a boundary point is polyhedral if and only if the corresponding direction cone is finitely generated.
Let us call such discrete point sets whose entire direction cones are finitely generated \emph{locally finitely generated} (cf. Definition \ref{Def:lfg}). Then it is possible to prove the following statement (cf. Theorem \ref{Thm4}):
\emph{All Voronoi cells of a discrete point set $\mathcal{P}$ are polyhedral if and only if $\mathcal{P}$ is locally finitely generated}. 

This paper which is a part of a Ph.D. thesis \cite{Voigt08} is organized as follows: 
In Section 2 we provide some well known definitions and facts from convex geometry, which will be applied below.
We also present an example of a non-polyhedral Voronoi cell of an infinite discrete point set.
In Section 3 we start our analysis of the shape of the Voronoi cells.
We introduce the terms of \emph{inner points} and \emph{boundary points} of a discrete point set and prove that the Voronoi diagram contains only polytopal cells if and only if all points of the discrete point set are inner points.
Subsequently, we investigate the behavior of the Voronoi cells of boundary points in Section 4. 
We define the \emph{direction cone} of a point of a discrete point set and show that the Voronoi cell of a boundary point is polyhedral if and only if its direction cone is finitely generated.
Overall, this yields our desired characterization.
Section 5 provides a short conclusion.

\section{Basic Principles and Settings}
Let $\mathbb{E}^n$ be the Euclidean space of dimension $n$ and denote the Euclidean norm by $\| . \|$ and the standard scalar product by $\langle .,. \rangle$.
We define the \emph{interior} of a set $M \subset \mathbb{E}^n$ by $M^{\circ} := \bigcup_{A \subset M, A \text{ open}}A$. The \emph{closure} of $M$ is given by $\overline{M} := \bigcap_{A \supset M, A \text{ closed}}A$. 
By defining the \emph{complement} of $M$ as $M^c := \mathbb{E}^n \smallsetminus M$, we get the \emph{boundary} of $M$ as $\partial M := \overline{M} \cap \overline{M^c}$. 
Furthermore, the \emph{affine hull} of $M$ is denoted by $\aff(M) := \bigcap_{A \supset M, A \text{ affine}}A$, whereat a set $M \subset \mathbb{E}^n$ is called \emph{affine} if $\lambda x + (1 - \lambda)y \in M$ for all $x,y \in M$ and $\lambda \in \mathbb{R}$.
We call a point set $\mathcal{P} \subset \mathbb{E}^n$ \emph{discrete} if it has only finitely many points in any bounded set of $\mathbb{E}^n$. 

For each point $p$ of a discrete point set $\mathcal{P} \subset \mathbb{E}^n$, we can define the \emph{Voronoi cell} $V(p)$ of $p$ relative to $\mathcal{P}$ as 
$$V(p) := \{ x \in \mathbb{E}^n \mid \| x - p \| \leq \| x - q \| \; \text{ for all } q \in \mathcal{P} \}.$$
The associated tessellation of $\mathbb{E}^n$ by all Voronoi cells of the set $\mathcal{P}$ is called the \emph{Voronoi diagram} of $\mathcal{P}$.
Each cell $V(p)$ ist the \emph{nearest neighbor region} of $p$ in $\mathbb{E}^n$.
An equivalent description of $V(p)$ is the representation of $V(p)$ as intersection of closed half spaces,
$$V(p) = \bigcap_{q \in \mathcal{P} \smallsetminus \{ p \}} H_p^-(q),$$
where the half spaces $H_p^-(q)$ are defined as $H_p^-(q) := H^-(q-p, \frac{1}{2}\| q-p \|^2)+p$ with $H^-(y, \frac{1}{2}\| y \|^2) := \{ x \in \mathbb{E}^n \mid \langle x,y \rangle \leq \frac{1}{2}\| y \|^2 \}$. 
Thus, the hyperplane $\Hyp_p(q)$ that bounds the half space $H_p^-(q)$ is orthogonal to the vector $q-p$ dividing the segment $[p,q]$ into equal parts.
Obviously, we do not always need all the half spaces $H_p^-(q)$ to represent $V(p)$ as an intersection.
Therefore, we call a point $p' \in \mathcal{P}$ (and a half space $H_p^-(p')$, respectively) \emph{Voronoi relevant} for $p$ if $V(p) \subsetneq \bigcap_{q \in \mathcal{P} \smallsetminus \{ p,p' \}} H_p^-(q)$.

As $V(p)$ can be described as an intersection of closed convex sets, the Voronoi cell itself is closed and convex.
Thereby, we call a set $K \subset \mathbb{E}^n$ \emph{convex}, if $(1-\lambda)x+\lambda y \in K$ for all $x,y \in K$ and $0 \leq \lambda \leq 1$.
If a convex set $K \subset \mathbb{E}^n$ can be represented as the intersection of finitely many closed half spaces, it is called a \emph{polyhedron}, and a bounded polyhedron is called a \emph{polytope}.
Hence, all Voronoi cells of a finite point set are polyhedra.

A point set $\mathcal{P} \subset \mathbb{E}^n$ is called a \emph{Delaunay set} (or $(r,R)$-\emph{system}), if there exist scalars $R > r > 0$ such that each ball of radius $r$ contains at most one point of $\mathcal{P}$, and every ball of radius $R$ contains at least one point of $\mathcal{P}$.
\textsc{Engel} showed in  Section 2 of \cite{Engel86} that each Voronoi cell of a Delaunay set is a polytope, and in particular, it is a polyhedron.

\begin{figure}
	\begin{center}
		\includegraphics[scale=0.75]{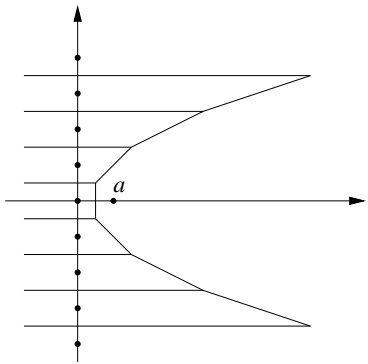}
	\end{center}
	\caption{Example \ref{Ex:non-polyhedral-cell}, non-polyhedral Voronoi cell}
	\label{fig:bsp3_2_1}
\end{figure}

If a closed convex set $K \subset \mathbb{E}^n$ has the property that each intersection of $K$ with a polytope is again a polytope, we call $K$ a \emph{generalized polyhedron}.
It is already known (cf. Chapter 32 of \cite{Grub07}) that all Voronoi cells of an arbitrary discrete point set are generalized polyhedra.
Noteworthy, an arbitrary generalized polyhedron is not always polyhedral, as we will see in the following example, where a discrete point set possessing a non-polyhedral Voronoi cell is given:

\begin{Ex}\label{Ex:non-polyhedral-cell}
\emph{
We set $\mathcal{P}_1 := \{ (0,z) \mid z \in \mathbb{Z} \} \cup \{ (1,0) \}$ and denote $a := (1,0)$. 
Then it is obvious that all half spaces $H_a^-(q), \; q \in \mathcal{P}_1 \smallsetminus \{ a \}$ and thus all points $q \in \mathcal{P}_1 \smallsetminus \{ a \}$ are Voronoi relevant for $a$. Thus $V(a)$ is not polyhedral, cf. Figure \ref{fig:bsp3_2_1}.
}
\end{Ex}

We now want to characterize the discrete point sets, for which all Voronoi cells are either polytopes or polyhedra.
For this task, we will need some well known definitions and facts from convex geometry which we will compile in the following. 
For the proofs, we refer to the literature on convex geometry like \cite{Barv02} and \cite{Grub07}. \\

For any set $M \subset \mathbb{E}^n$ we can construct the \emph{convex hull} $\conv(M)$ as the set, which consists of all convex combinations of points of $M$ that is
$$\conv(M) := \Bigl\{ \sum_{i=1}^k \lambda_i x_i \mid x_i \in M, \; \lambda_i \geq 0, \; \sum_{i=1}^k \lambda_i = 1, \; k \in \mathbb{N} \Bigr\}.$$
Furthermore, by \textsc{Carath\'eodory}'s Theorem we know that each point $x$ of $\conv(M)$ can be represented as a convex combination of at most $n+1$ points of $M$ (cf. Theorem I.2.3 in \cite{Barv02}).

\begin{Thm}\label{Thm:Caratheodory}
\textbf{(Carath\'eodory)}
For $M \subset \mathbb{E}^n$ and $x \in \conv(M)$ there exist $a_1, \dots, a_k \in M$ such that $x = \sum_{i=1}^k \lambda_i a_i$ with $\lambda_1, \dots, \lambda_k \geq 0, \; \sum_{i=1}^k \lambda_i = 1$ and $k \leq n+1$. 
\end{Thm}

A similar statement holds for interior points of the convex hull; \textsc{Steinitz} showed (cf. Problem No. 3 to I.2.3 in \cite{Barv02}):

\begin{Thm}\label{Thm:Steinitz}
\textbf{(Steinitz)} 
For $M \subset \mathbb{E}^n$ and $x \in \conv(M)^{\circ}$ there exists a subset $M' \subset M$ with $|M'| \leq 2n$ such that $x \in \conv(M')^{\circ}$.
\end{Thm}

Beside the convex hull, we can also construct a \emph{conic hull} of any given set $M \subset \mathbb{E}^n$,
$$\cone{M} := \Bigl\{ \sum_{i=1}^k \lambda_i x_i \mid x_i \in M, \; \lambda_i \geq 0, \; k \in \mathbb{N} \Bigr\},$$
which is a convex cone with apex $0$. 

In order to characterize convex sets, it is somehow convenient to look for the greatest cone which is contained in the considered convex set.
Therefor we define the \emph{recession cone} (also called \emph{characteristic cone}) of a convex set $K \subset \mathbb{E}^n$ as
$$C_c(K) := \{ x \in \mathbb{E}^n \mid x+y \in K \text{ for all } y \in K \}.$$
For each convex set $K$, its recession cone $C_c(K)$ is either a convex cone with apex $0$ or equals $\{ 0 \}$.
In the special case that the convex set $K$ is a polyhedron, $K = P = \bigcap_{i=1}^k \{ x \in \mathbb{E}^n \mid \langle x,y_i \rangle \geq \alpha_i \}$, we get $C_c(P) = \bigcap_{i=1}^k \{ x \in \mathbb{E}^n \mid \langle x,y_i \rangle \geq 0 \}$. 

Further investigation of convex sets deal with the determination of extreme points and rays.
For a convex set $K \subset \mathbb{E}^n$, we call $x \in K$ an \emph{extreme point} of $K$ if there exists no non-trivial convex combination $x = (1-\lambda)x_1 + \lambda x_2$ with points $x_1,x_2 \in K$ and $0 < \lambda < 1$.
The set of all extreme points of $K$ is denoted by $\Ext(K)$.
By a theorem of \textsc{Krein} and \textsc{Milman},
we know that each closed, bounded convex set is uniquely determined by its extreme points (cf. Theorem III.4.1 in \cite{Barv02}).

\begin{Thm}\label{Thm:Krein-Milman}
\textbf{(Krein-Milman)}
If $K \subset \mathbb{E}^n$ is a closed, bounded convex set, then $K = \conv(\Ext(K))$.
\end{Thm}

It is obvious that the determination of $\Ext(K)$ is no longer sufficient, if $K$ is unbounded. In this case, we have to comprise rays that belong to $K$. Thereby, a \emph{ray} with \emph{initial point} $x_0$ and \emph{direction} $v$ is defined as 
$$S(x_0,v) := \{ x \in \mathbb{E}^n \mid x = x_0 + \lambda v, \lambda \geq 0 \}.$$
We call a ray $S(x_0,v) \subset K$ of a convex set $K$ an \emph{extreme ray} of $K$, if $x_0$ is an extreme point of $K$ and $K \smallsetminus S(x_0,v)$ is still convex.
The set of all extreme rays of $K$ is denoted by $\RExt(K)$.
\textsc{Klee} proved that each line-free, closed convex set is generated by its extreme points and extreme rays (cf. Lemma II.16.3 in \cite{Barv02}).

\begin{Thm}\label{Thm:Klee}
\textbf{(Klee)}
Each closed convex set $K \subset \mathbb{E}^n$ that contains no line fulfills $K = \conv(\Ext(K) \cup \RExt(K))$.
\end{Thm}

Using this result, one can show that each closed convex set has a decomposition as Minkowski sum of an affine subspace together with a convex cone and a compact convex set.

Another approach to describe a convex set $K \subset \mathbb{E}^n$ utilizes the hyperplanes which are tangential to $K$. Therefore, we call a hyperplane $\Hyp$ \emph{supporting hyperplane} to $K$ if $K$ is completely contained in one of the half spaces bounded by $\Hyp$ and $\Hyp \cap K \neq \emptyset$. 
The intersection $F := \Hyp \cap K$ is called a \emph{face} of $K$ if $\Hyp$ is a supporting hyperplane of $K$.

It can be shown that it is possible to find for every convex set $K \subset \mathbb{E}^n$ a \emph{dual} convex set such that the inclusion-chains of the faces are reversed. One realization of such a dual set is the \emph{polar set} $K^*$, which we define as $$K^* := \{ y \in \mathbb{E}^n \mid \langle x,y \rangle \leq 1 \text{ for all } x \in K \}.$$ 
We subsume the most important properties of polar sets (cf. Chapter IV, Section 1 in \cite{Barv02}):

\begin{Thm}\label{Thm:polar-set}
Let $K \subset \mathbb{E}^n$ be a convex set.
\begin{enumerate}
\item[1. ] If $K$ is compact and the origin $0 \in K^{\circ}$ is an interior point, the same holds for $K^*$, and we have $K^{**} = K$.
\item[2. ] $K$ is bounded if and only if $0 \in (K^*)^{\circ}$.
\item[3. ] For all convex sets $K_1 \subset K_2 \subset \mathbb{E}^n$, it holds $K_2^* \subset K_1^*$.
\item[4. ] If $K = \conv(y_1, \dots, y_k)$ is a polytope, with $y_1, \dots, y_k \in \mathbb{E}^n$, we get $K^* = \{ x \in \mathbb{E}^n \mid \langle x,y_i \rangle \leq 1 \text{ for } i = 1, \dots, k \}$.
\end{enumerate}
\end{Thm}

If $C \subset \mathbb{E}^n$ is a convex cone with apex $0$, the polar cone is equal to
$$C^* = \{ y \in \mathbb{E}^n \mid \langle x,y \rangle \leq 0 \text{ for all } x \in C \}.$$

\section{Voronoi cells of inner points}
In the above Example \ref{Ex:non-polyhedral-cell}, we observe that the convex hull $\conv(\mathcal{P}_1)$ is not closed and that the non-polyhedral Voronoi cell $V(a)$ belongs to a boundary point of $\conv(\mathcal{P}_1)$.
We want to compare this with the cases of Delaunay sets, infinite discrete point sets in general quadratic position and finite sets:
If the point set is a Delaunay set with parameters $(r,R)$, the existence of the scalar $R$ yields that all points belong to the interior of the convex hull of the point set, which equals the whole space. And it is known that all Voronoi cells of a Delaunay set are polytopes (cf. Section 2 of \cite{Engel86}).
In case of an infinite discrete point set in $\mathbb{E}^n$ which is in general quadratic position, we get that all subsets of $n+1$ points are affinely independent. Hence its convex hull equals the whole space and each point is an interior point. We also know that all Voronoi cells are polytopes (cf. Chapter 5 of \cite{OBS00}).
If we consider a finite point set, we see that the Voronoi cells of points that are contained in the interior of its convex hull are again polytopes, while we only get polyhedral cells if the corresponding point lies on the boundary of the convex hull (cf. Property V2 of Chapter 2.3 in \cite{OBS00}).
In light of the above, we aim to prove the following theorem:
For each (infinite) discrete point set holds, that the Voronoi cell of a point is a polytope if and only if the point belongs to the interior of the convex hull.

To simplify the notation, we want to introduce the terms of \emph{inner points} and \emph{boundary points} of discrete sets:

\begin{Def}\label{Def:inner+boundary-points}
\emph{
Let $\mathcal{P} \subset \mathbb{E}^n$ be a discrete point set.
We call a point $p \in \mathcal{P}$ an \emph{inner point} of $\mathcal{P}$, if $p \in \conv(\mathcal{P})^{\circ}$.
The set of all inner points of $\mathcal{P}$ is denoted by $\mathcal{P}^{\circ} := \mathcal{P} \cap \conv(\mathcal{P})^{\circ}$.
A point $p' \in \mathcal{P}$ is called a \emph{boundary point} of $\mathcal{P}$, if $p' \in \partial \conv(\mathcal{P})$.
We denote the set of all boundary points by $\partial \mathcal{P} := \mathcal{P} \cap \partial \conv(\mathcal{P})$.
}
\end{Def}

\begin{figure}
	\begin{center}
		\includegraphics[scale=0.75]{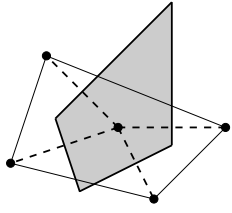}
	\end{center}
	\caption{Voronoi cell of an inner point}
	\label{fig:vor_art_thm1}
\end{figure}

\begin{Thm}\label{Thm1}
Let $\mathcal{P} \subset \mathbb{E}^n$ be a discrete point set. Then the Voronoi cell $V(p)$ of a point $p \in \mathcal{P}$ is a polytope if and only if $p \in \mathcal{P}^{\circ}$ (cf. Figure \ref{fig:vor_art_thm1}).
\end{Thm}

\textbf{Proof.} 
Without loss of generality, we assume that the point $p \in \mathcal{P}$ lies in the origin $p = 0$. (Otherwise we consider the translated point set $\mathcal{P}-p = \{ q-p \mid q \in \mathcal{P} \}$.)

Now, let $p \in \mathcal{P}^{\circ}$ be an inner point of the discrete point set. 
Then Steinitz's Theorem \ref{Thm:Steinitz} yields that there exist finitely many points $p_1, \dots, p_k$ in $\mathcal{P}$ with $k \leq 2n$ such that $p$ is an interior point of $\conv(p_1, \dots, p_k)$.
Let $C := \conv(p_1, \dots, p_k)^*$ be the polar set of $\conv(p_1, \dots, p_k)$.
Then by applying Theorem \ref{Thm:polar-set}, we obtain $C = \{ x \in \mathbb{E}^n \mid \langle x,p_i \rangle \leq 1, \;i = 1, \dots, k \}$.
Since $\conv(p_1, \dots, p_k)$ is compact and $p$ is an interior point of $\conv(p_1, \dots, p_k)$, the polar set $C$ is compact and contains $p$ as an interior point.
Setting $m := \max \{ \frac{1}{2}\|p_i\|^2 \mid i=1,\dots,k \}$, we obtain $\bigcap_{i=1}^k H_p^-(p_i) \subset m \cdot C = \{ mx \mid x \in C \}$.
As the Voronoi cell $V(p) = \bigcap_{q \in \mathcal{P} \smallsetminus \{ p \}} H_p^-(q)$, we obviously get $V(p) \subset m \cdot C$. Thus $V(p)$ is contained in a polytope.
Owing to the fact that each Voronoi cell is a generalized polyhedron, we obtain that $V(p) = V(p) \cap m \cdot C$ is a polytope.

To show the converse, let $V(p) = \bigcap_{i=1}^k H_p^-(p_i)$ be a polytope.
Then, Theorem \ref{Thm:polar-set} yields that $V(p)^*$ is also a polytope with the origin $p=0$ as an interior point.
In addition, we know that $V(p)^*$ is of the form $V(p)^* = \conv(q_1, \dots, q_k)$ for suitable $q_i = \lambda_i p_i, \, \lambda_i > 0$ for $i = 1, \dots, k$.
Thus it follows directly that $p$ is an interior point of $\conv(p_1, \dots, p_k)$ and in particular $0=p \in \mathcal{P}^{\circ}$.
\hfill $\Box$ \\

We have just characterized the discrete point sets, for which the associated Voronoi diagram contains only polytopal cells, as exactly those point sets that consist only of inner points.
Furthermore, as the convex hull of a Delaunay set or an infinite discrete point set in general quadratic position equals the whole space $\mathbb{E}^n$, each point of these sets is an inner point, and we get as a special case of Theorem \ref{Thm1} that all Voronoi cells are polytopes.

\begin{Cor}\label{Cor:1}
Let $\mathcal{P} \subset \mathbb{E}^n$ be a discrete point set. Then it holds:
\begin{enumerate}
\item[(1) ] All Voronoi cells of $\mathcal{P}$ are polytopes if and only if $\mathcal{P} = \mathcal{P}^{\circ}$.
\item[(2) ] If $\conv(\mathcal{P}) = \mathbb{E}^n$, then all Voronoi cells of $\mathcal{P}$ are polytopes.
\item[(3) ] If the Voronoi cell $V(p)$ of a point $p \in \mathcal{P}$ is non-polyhedral, then $p$ belongs to $\partial \mathcal{P}$.

\end{enumerate}
\hfill $\Box$
\end{Cor}

Interestingly, for a discrete point set $\mathcal{P}$ the property $\mathcal{P} = \mathcal{P}^{\circ}$ is not equivalent to $\conv(\mathcal{P}) = \mathbb{E}^n$. This can be illustrated by the following example:

\begin{figure}
	\begin{center}
		\includegraphics[scale=0.75]{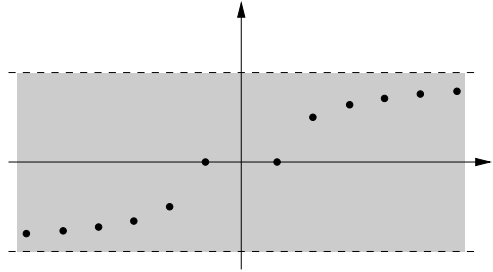}
	\end{center}
	\caption{Example \ref{Ex:only-inner-points}, convex hull}
	\label{fig:folg3_2_10plus}
\end{figure}

\begin{Ex}\label{Ex:only-inner-points}
\emph{
Consider the point set $\mathcal{P}_2 := \{ (n,1-\frac{1}{n}) \mid n \in \mathbb{N} \} \cup \{ (-n,-1+\frac{1}{n}) \mid n \in \mathbb{N} \} \subset \mathbb{E}^2$, where each point of $\mathcal{P}_2$ is an inner point, cf. Figure \ref{fig:folg3_2_10plus}. It follows from Corollary \ref{Cor:1} that all Voronoi cells are polytopes, but the convex hull of $\mathcal{P}_2$ does not equal $\mathbb{E}^2$.
}
\end{Ex}

\begin{figure}
	\begin{center}
		\includegraphics[scale=0.75]{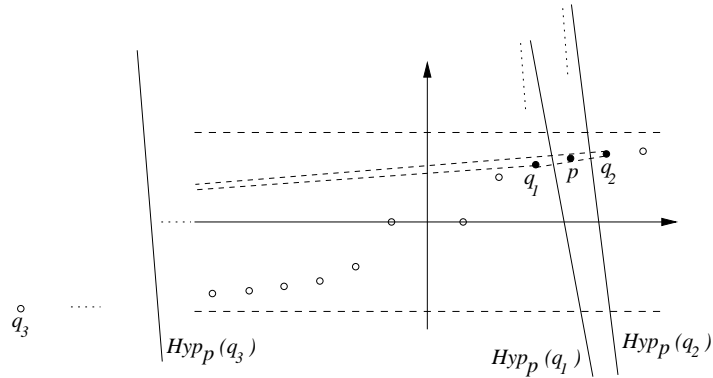}
	\end{center}
	\caption{Example \ref{Ex:only-inner-points}, Voronoi cell construction}
	\label{fig:folg3_2_10}
\end{figure}

If we now want to sketch the Voronoi diagram of the set $\mathcal{P}_2$ of Example \ref{Ex:only-inner-points} in order to see all the polytopal Voronoi cells, we recognize that the cells exhibit a very "strange" shape. That is, they are very long and thin, and it is hard to believe that all of them are polytopes.
Figure \ref{fig:folg3_2_10} illustrates that this is nevertheless true. 
Let us choose an arbitrary point $p := (x,y) \in \mathcal{P}_2$, without loss of generality $x,y > 0$.
To construct the Voronoi cell $V(p)$, we first consider the hyperplanes $\Hyp_p(q_1), \Hyp_p(q_2)$ corresponding to the two closest points $q_1, q_2 \in \mathcal{P}_2$ to $p$.
Since these hyperplanes are not parallel, they intersect in a point $(x_1,-y_1)$ with $y_1 \gg 0$.
Furthermore, we find another point $q_3 \in \mathcal{P}_2$ far away from $p$ whose corresponding hyperplane $\Hyp_p(q_3)$ intersects $\Hyp_p(q_1), \Hyp_p(q_2)$ in points $(x_2,y_2), \; (x_3,y_3)$ with $y_2,y_3 \gg 0$. Therefore, the Voronoi cell $V(p)$ is bounded. Hence, we have also geometrically confirmed that the Voronoi cells of $\mathcal{P}_2$ are indeed polytopes. \\

\section{Voronoi cells of boundary points}

To construct the Voronoi cell $V(p)$ of a point $p$ of a discrete point set $\mathcal{P}$, we have to consider the half spaces, which are defined by the directions $q-p$ for $q \in \mathcal{P} \smallsetminus \{ p \}$. 
It is therefore convenient to investigate the cone which is generated by the directions emanating from $p$. We want to call this cone the \emph{direction cone} of $p$.

\begin{Def}\label{C(p)}
\emph{
Let $\mathcal{P} \subset \mathbb{E}^n$ be discrete and $p \in \mathcal{P}$.
Then we define the \emph{direction cone} of $p$ by
$$ C(p) := \cone(\mathcal{P}-p) = \Bigl\{ \sum_{i=1}^m \lambda_i (p_i-p) \mid p_i \in \mathcal{P}, \; \lambda_i \geq 0, \; m \in \mathbb{N} \Bigr\}. $$ 
}
\end{Def}

One easily derives that the direction cone $C(p)$ of a point $p$ is either a convex cone with apex $0$ or the whole space $\mathbb{E}^n$.
More precisely, it is possible to prove that the conic hull $C(p)$ of a point $p$ equals the whole space $\mathbb{E}^n$ if and only if $p$ is an inner point.
It is also clear that, for $p \in \partial \mathcal{P}$, the directions which define $C(p)$ are perpendicular to the hyperplanes of the half spaces which define the Voronoi cell $V(p)$. 
Thus, there is some kind of "duality" between $C(p)$ and $V(p)$.
To specify this duality, we will prove that the recession cone of a polyhedral Voronoi cell $V(p)$ of a boundary point $p \in \partial \mathcal{P}$ equals the polar cone of the direction cone $C(p)$, that is, $C(p)^* = C_c(V(p))$.
This yields the following relation: If the Voronoi cell of $p$ is finitely generated (i.e. if $V(p)$ is polyhedral), then the same holds for the direction cone $C(p)$.

\begin{Thm}\label{Thm2}
Let $\mathcal{P} \subset \mathbb{E}^n$ be a discrete point set and $p \in \partial \mathcal{P}$ a boundary point. 
Then $C(p)$ is finitely generated, if $V(p)$ is polyhedral.
\end{Thm}

\textbf{Proof.}
As $V(p)$ is polyhedral, there exist finitely many points $p_1, \dots, p_k$ in $\mathcal{P}$ such that
$$ V(p) = \bigcap_{q \in \mathcal{P} \smallsetminus \{ p \}} H_p^-(q) = \bigcap_{i=1}^k H_p^-(p_i) = \bigcap_{i=1}^k H^-(p_i-p,\frac{1}{2}\| p_i-p \|^2) + p. $$
Thus one obtains
$$ C_c(V(p)) = \bigcap_{i=1}^k H^-(p_i-p,0) = \bigcap_{q \in \mathcal{P} \smallsetminus \{ p \}} H^-(q-p,0), $$ 
and moreover \\
$$\begin{array}{lll} 
\cone(p_1-p, \dots, p_k-p) & = (\bigcap_{i=1}^k H^-(p_i-p,0))^* & \\
 & = (C_c(V(p)))^* & \\
 & = (\bigcap_{q \in \mathcal{P} \smallsetminus \{ p \}} H^-(q-p,0))^* & \\
 & = \cone(\mathcal{P}-p) & = C(p). 
\end{array}$$ \\
Hence, it has been shown that the direction cone is finitely generated.
\hfill $\Box$ \\

A finitely generated cone is obviously closed, and hence, we can conclude that the Voronoi cell $V(p)$ of a boundary point $p \in \partial \mathcal{P}$ is non-polyhedral, if the direction cone $C(p)$ is not closed. 
The Voronoi cell $V(a)$ of Example \ref{Ex:non-polyhedral-cell} provides an example for this property.

\begin{Thm}\label{Thm3}
Let $\mathcal{P} \subset \mathbb{E}^n$ be a discrete point set and $p \in \partial \mathcal{P}$ a boundary point. 
Then $V(p)$ is polyhedral if $C(p)$ is finitely generated.
\end{Thm}

\begin{figure}
	\begin{center}
		\includegraphics[scale=0.75]{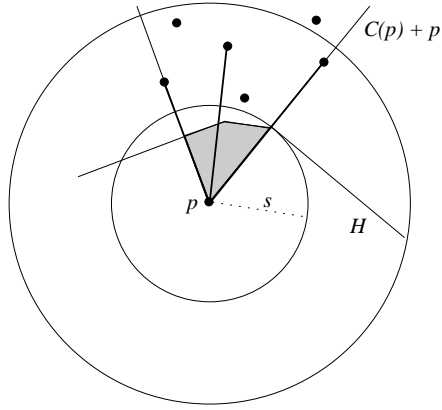}
	\end{center}
	\caption{Proof of Theorem \ref{Thm3}, construction}
	\label{fig:satz3_4_15}
\end{figure}

\textbf{Proof.}
As $C(p)$ is finitely generated, we find $p_1, \dots, p_k \in \mathcal{P}$ with 
$$ C(p) = \cone(\mathcal{P} - p) = \cone(p_1-p, \dots, p_k-p). $$ 
If we set 
$$ H := \bigcap_{i=1}^k H_p^-(p_i), $$ 
we obviously obtain $V(p) \subset H$ and 
$$ C_c(H) = \bigcap_{i=1}^k H^-(p_i-p,0) = C(p)^*. $$
Since $p$ is a boundary point, we have $C(p) \neq \mathbb{E}^n$. 
Thus, $C(p)^* = C_c(H) \neq \{ 0 \}$, so that $C(p)$ is not bounded.
Note that a Voronoi relevant point must be in $C(p)+p$ and its associated hyperplane must have a non-empty intersection with the interior of $H$. It is therefore convenient to consider the set 
$$ H \cap (C(p)+p) $$
(cf. Figure \ref{fig:satz3_4_15}). 
Owing to the fact that $C_c(H) \cap C_c(C(p)+p) = C(p)^* \cap C(p) = \{ 0 \}$, it can easily be seen that $H \cap (C(p)+p)$ is bounded. Since $H \cap (C(p)+p)$ is the intersection of two closed sets, $H \cap (C(p)+p)$ is as well closed and thus compact.
Hence, it exists a maximal distance 
$$ s' := \max \{ \|x-p\| \mid x \in \partial(H \cap C(p)+p) \} < \infty $$ 
between the point $p$ and the boundary of $H \cap (C(p)+p)$.
Since $H$ is a polyhedron, the number of extreme points $|\Ext(H)| < \infty$ of $H$ is finite.
Therefore, we can also define 
$$ s'' := \max \{ \|x-p\| \mid x \in \Ext(H) \} < \infty, $$
and we set 
$$ s := \max \{ s',s'' \}. $$
This allows us to show that all points, which are relevant for $p$, belong to the ball $B(p,2s) := \{ x \in \mathbb{E}^n \mid \|x-p\| \leq 2s \}$.
We prove this by contradiction: 

We assume that we have a point $q \in \mathcal{P} \subset C(p)+p$ which is relevant for $p$, but $q \not\in B(p,2s)$.
As the bounding hyperplane $\Hyp_p(q)$ of the associated half space $H_p^-(q)$ is tangent to the ball $B(p,\frac{1}{2}\|q-p\|)$ and $\frac{1}{2}\|q-p\| > s$, we get $\Hyp_p(q) \cap B(p,s) = \emptyset$, and thus, $\Hyp_p(q) \cap (H \cap (C(p)+p)) = \emptyset$.
Since $\Ext(H) \subset (B(p,s) \cap H) \subset (H \cap H_p^-(q))$ and $\frac{1}{2}\|q-p\| > s$, we get $\Ext(H) \subset \Ext(H \cap H_p^-(q))$.
Without loss of generality, we assume that $H$ is line-free (otherwise we restrict our attention to $H \cap \aff(C(p)+p)$, which is line-free).
With $H \cap H_p^-(q) \subsetneq H$, Theorem \ref{Thm:Klee} of Klee yields the inclusion
$$ \begin{array}{ll}
 & \conv(\Ext(H) \cup \RExt(H \cap H_p^-(q))) \\
\subseteq & \conv(\Ext(H \cap H_p^-(q)) \cup \RExt(H \cap H_p^-(q))) \quad = \quad H \cap H_p^-(q) \\ 
\subsetneq & H \quad = \quad \conv(\Ext(H) \cup \RExt(H)).
\end{array} $$
Consequently $\RExt(H) \smallsetminus \RExt(H \cap H_p^-(q)) \neq \emptyset$, and moreover we get $C_c(H) \smallsetminus C_c(H \cap H_p^-(q)) \neq \emptyset$.
As a result, we have $C_c(H \cap H_p^-(q)) \subsetneq C_c(H)$, what yields $C_c(H)^* \subsetneq C_c(H \cap H_p^-(q))^*$ by Theorem \ref{Thm:polar-set}.
Since
$$ C_c(H)^* = C(p) = \cone(p_1-p, \dots, p_k-p) $$
and  
$$ C_c(H \cap H_p^-(q))^* = \cone(p_1-p, \dots, p_k-p, q-p), $$ 
we conclude that $q \not\in C(p)+p$, and thus $q \not\in \mathcal{P}$; this is a contradiction.

Hence we have shown that $B(p,2s)$ contains all points which are relevant for $p$.
Since $\mathcal{P}$ is discrete, we have $|\mathcal{P} \cap B(p,2s)| < \infty$, so that
there exist only finitely many relevant points for $p$. Therefore, $V(p)$ is a polyhedron.
\hfill $\Box$ \\

Since the direction cone of an inner point is obviously finitely generated, we can extend the statements of Theorem \ref{Thm2} and \ref{Thm3} to all points of a discrete point set.
We can therefore conclude that all Voronoi cells of a discrete point set are polyhedral if and only if the direction cones of all points are finitely generated.
To emphasize this property we want to call such sets \emph{locally finitely generated}.

\begin{Def}\label{Def:lfg}
\emph{
Let $\mathcal{P} \subset \mathbb{E}^n$ be a discrete point set. We call $\mathcal{P}$ \emph{locally finitely generated} if the direction cones $C(p)$ of all points $p \in \mathcal{P}$ are finitely generated.
}
\end{Def}


\begin{Thm}\label{Thm4}
All Voronoi cells of a discrete point set $\mathcal{P}$ are polyhedral if and only if $\mathcal{P}$ is locally finitely generated.
\hfill $\Box$
\end{Thm}

Hence, we have achieved our aim and characterized the discrete point sets whose Voronoi diagram contains only polyhedral cells.
The formerly investigated cases of finite sets, Delaunay sets and infinte discrete point sets in general quadratic position fit in as special cases of the above theorem.

\section{Conclusion}
In the past, Voronoi diagrams were usually studied for finite sets or infinite sets with a certain structure, that is, the formerly considered infinite sets were either Delaunay sets or discrete point sets in general quadratic position.
One knows that in all these cases all Voronoi cells are polytopes or polyhedra.
By stating Example \ref{Ex:non-polyhedral-cell} we saw that this has not to be true for an arbitrary infinite discrete point set; there, the corresponding Voronoi diagram contained a non-polyhedral cell. 
Hence, we aimed to characterize the discrete point sets for which all Voronoi cells are polytopes or polyhedra.
We investigated which property of the discrete point set yields polyhedral or polytopal Voronoi cells and found the desired characterizations, see the Theorems \ref{Thm1} and \ref{Thm4}.
Finally, we identified all formerly considered cases as special cases of our theory. \\

\medskip
\noindent
\textbf{Acknowledgments.}
I want to thank Rudolf Scharlau for supervising my Ph.D. thesis \cite{Voigt08} and especially J\"urgen Eckhoff and Frank Vallentin for many useful tips and hints regarding the writing of this paper.

\end{document}